\documentclass[9pt]{article}
\usepackage{mathptmx}       
\usepackage{helvet}         
\usepackage{courier}        
\usepackage{type1cm}        
\usepackage{makeidx}         
\usepackage{graphicx}        
\usepackage{multicol}        
\usepackage[bottom]{footmisc}

\usepackage{geometry}                
\usepackage{graphicx}
\usepackage{amssymb}
\usepackage{amsmath}
\usepackage{epstopdf}
\usepackage{stmaryrd}
\DeclareMathAlphabet{\mathcal}{OMS}{cmsy}{m}{n}

\usepackage{accsupp} 

\newcommand\iprod[1]{\left( #1\right)} 				
\newcommand\inorm[1]{\left |\left| #1\right|\right|}		
\newcommand\iprodN[1]{\left( #1\right)_{N}}			
\newcommand\spacevec[1]{\vec{ #1}}					
\newcommand\contravec[1]{\tilde{ #1}}					
\newcommand\statevec[1]{\mathbf #1}					
\newcommand\mmatrix[1]{\mathcal{#1}}				
\newcommand{\jump}[1]{\left\llbracket#1\right\rrbracket}   
\newcommand{\partialderivative}[2]{\frac{\partial #1}{\partial #2}}  

\title{A Polynomial Spectral Calculus for Analysis of DG Spectral Element Methods}
\author{David A. Kopriva\\Department of Mathematics\\ The Florida State University\\ Tallahassee, FL 32306 USA\\ email:kopriva@math.fsu.edu}

\begin{document}
\maketitle
\abstract{We introduce a polynomial spectral calculus that follows from the summation by parts property of the Legendre-Gauss-Lobatto quadrature. We use the calculus to simplify the analysis of two multidimensional discontinuous Galerkin spectral element approximations.}

\section{Introduction}
The discontinuous Galerkin Spectral Element Method (DGSEM) introduced by Black \cite{Black1999},\cite{black2000} has the desired properties of spectral accuracy, geometric flexibility, and excellent phase and dissipation properties \cite{Stanescuetal2001},\cite{gassner2011}. Spectral accuracy comes from the use of high order polynomial approximations to the solutions and fluxes, and high order Gauss quadratures for the inner products, e.g. \cite{Kopriva2017}. Geometric flexibility comes from the multi-element subdivision of the domain. The DGSEM is now developed to the point of being efficient for large scale engineering level computations, e.g. \cite{Altmann2012},\cite{FLD:FLD3943},\cite{Frank2016132}, among others.

Robustness, however, has been an issue with the DGSEM at high order. It usually works, but it can go unstable even when the solutions are smooth. For nonlinear problems, this is probably not surprising. Examples are demonstrated in the computation of the Taylor-Green vortex problem, where instability at high orders is seen \cite{Gassner:2013qf}. But instability arises even in linear problems when the coefficients are variable, which can come from inherent variability \cite{Beck:2013sf} or from variability introduced by curved elements \cite{Gassner:2013ij}. The instability, we will show, comes from aliasing errors associated with the products of polynomials and  insufficient Gauss quadrature precision.

Robust (provably stable) versions of the DGSEM that start from a split form of the partial differential equation (PDE) have recently been developed for linear hyperbolic systems for static \cite{Gassner:2013ij} and moving domains \cite{Kopriva:2016il}. In addition to stability, the approximations match the additional conservative and constant state preserving properties of the PDE \cite{Kopriva2016274}. The approach is applicable to nonlinear problems, where, depending on the equations and split form, the methods are energy or entropy stable \cite{gassner_skew_burgers},\cite{Gassner2016291},\cite{Gassner:2016ye}.

In this paper, we introduce a polynomial spectral calculus that allows us to mirror the continuous PDE analysis to show stability of Black's and the split-form approximations with a simple, compact notation applicable to any number of space dimensions. For the split form method, we also show how to use the calculus to demonstrate conservation and constant state preservation. The key starting point of the calculus is the summation by parts property satisfied by the Gauss-Lobatto quadrature \cite{gassner2010}, which allows us to write discrete versions of the Gauss law and its variants. Those discrete Gauss laws, in turn, allow us to write algebraically equivalent forms of the approximations, with which we can easily analyze their properties.

\section{Linear Hyperbolic Problems on Bounded Domains}

As examples of the use of the discrete calculus, we will analyze two discontinuous Galerkin spectral element approximations to the linear system of conservation laws 
\begin{equation}
{{\mathbf{u}}_t} + \nabla  \cdot \spacevec {\mathbf{f}}  = 0,
\label{eq:ConservativeSystem}
\end{equation}
where $\mathbf{u}\left(\spacevec x,t\right) = [u_{1}\;u_{2}\;\ldots \; u_{p}]^{T}$ is the state vector and 
\begin{equation}
\spacevec {\mathbf{f}} \left(\mathbf u\right) = \sum\limits_{m = 1}^3 {{\mathcal{A}^{(m)}}\left( {\spacevec x} \right){\mathbf{u}}{{\hat x}_m}}  \equiv \spacevec {\mathcal{A}} {\mathbf{u}}
\end{equation}
is the linear flux space-state vector. For simplicity we will assume that the system has been symmetrized and is hyperbolic so that
\begin{equation}
\mmatrix{A}^{(m)} = \left(\mmatrix{A}^{(m)}\right)^{T}\quad\text{and}\quad \sum\limits_{m = 1}^3 {{\alpha _m}{\mmatrix{A}^{(m)}}}  = \mmatrix{R}\left(\spacevec\alpha\right)\Lambda\left(\spacevec\alpha\right) {\mmatrix{R}^{ - 1}}\left(\spacevec\alpha\right)
\label{eq:mDHyperbolicConditions_DAK}
\end{equation}
for any $\left\|\spacevec \alpha\right\|^{2}_{2} =\sum\limits_{m = 1}^3 {\alpha _m^2}  \ne 0$ and some real diagonal matrix $\Lambda$. We will also assume that the matrices $\mathcal{A}^{(m)}$ have bounded derivatives in the sense that 
\begin{equation}{\left\| {\nabla  \cdot \spacevec{ \mathcal{A}} } \right\|_2}<\infty,\end{equation}
where $\left\|\cdot\right\|_2$ is the matrix 2-norm. Additional constraints on the coefficient matrices need to be added later to ensure that the derivatives of their interpolants converge in the maximum norm.
 The product rule applied to (\ref{eq:ConservativeSystem} leads to the nonconservative form of the system
\begin{equation}
\mathbf{u}_{t}+\left(\nabla\cdot\spacevec{\mathcal{A}}\right)\mathbf{u} + \spacevec{\mathcal{A}}\cdot\nabla\mathbf{u}=0.
\label{eq:NonConForm}
\end{equation}

With appropriate initial and characteristic boundary conditions on a bounded domain $\Omega\in \mathbb R^{3}$ the problem is (i) well posed, (ii) conservative, and, under conditions on $\spacevec{\mathcal{A}}$, (iii) preserves a constant state. These properties are demonstrated from a weak form of the average of the conservative, (\ref{eq:ConservativeSystem}), and nonconservative, (\ref{eq:NonConForm}), forms of the equation, the so-called ``split-form''. To write the weak form, we define the $\mathbb{L}^{2}$ inner product and norm
\begin{equation}
\left( {{\mathbf{u}},{\mathbf{v}}} \right) = \int_\Omega  {{{\mathbf{u}}^T}{\mathbf{v}}dxdydz} ,\quad \left\| {\mathbf{u}} \right\| = \sqrt {\left( {{\mathbf{u}},{\mathbf{u}}} \right)}.
\end{equation}
Then for any state vector $\boldsymbol{\phi} \in \mathbb{L}^{2}\left(\Omega\right)$,
\begin{equation}
\left( {{{\mathbf{u}}_t},\boldsymbol\phi } \right) + \frac{1}{2} \left( {\nabla  \cdot \spacevec {\mathbf{f}} ,\boldsymbol\phi } \right) + \frac{1}{2}\left\{ {\left( {\left( {\nabla  \cdot \spacevec{ \mathcal{A}} } \right){\mathbf{u}},\boldsymbol\phi } \right) + \left( {\spacevec{ \mathcal{A}}  \cdot \nabla {\mathbf{u}},\boldsymbol\phi } \right)} \right\} = 0.
\label{eq:WeakForm1}
\end{equation}

From vector calculus, we have the extended Gauss law, 
\begin{equation}\int_\Omega  {{\statevec u^T}\nabla  \cdot \spacevec{\statevec f} dxdydz}  = \int_{\partial \Omega } {{\statevec u^T}\spacevec{\statevec f} \cdot \hat ndS}  - \int_\Omega  {{{\left( {\nabla \statevec u} \right)}^T} \cdot \spacevec{\statevec f} dxdydz}, \label{eq:ExtendGauss}\end{equation}
where $\hat n$ is the outward unit normal. We write (\ref{eq:ExtendGauss}) in inner product form as
\begin{equation}\left( {\statevec u,\nabla  \cdot \spacevec{\statevec f}} \right) = \int_{\partial \Omega } {{\statevec u^T}\spacevec{\statevec f} \cdot \hat ndS}  - \left( {\nabla \statevec u,\spacevec{ \statevec f}} \right).\end{equation}

We can apply the extended Gauss law to the inner products in the braces in (\ref{eq:WeakForm1}) and use the fact that $\spacevec{\mathcal{A}}$ is symmetric to get an equivalent form that separates the boundary and volume contributions
\begin{equation}
\left( {{{\mathbf{u}}_t},\boldsymbol\phi } \right) + \int_{\partial \Omega } {\spacevec {\mathbf{f}}  \cdot \hat n\boldsymbol\phi dS}  - \frac{1}{2}\left( {\spacevec {\mathbf{f}} ,\nabla \boldsymbol\phi } \right) + \frac{1}{2}\left\{ {\left( {\left( {\nabla  \cdot \spacevec{ \mathcal{A}} } \right){\mathbf{u}},\phi } \right) - \left( {{\mathbf{u}},\nabla  \cdot \left( {\spacevec{ \mathcal{A}} \boldsymbol\phi } \right)} \right)} \right\} = 0.
\label{eq:WeakForm2}
\end{equation}

Constant state preservation, conservation and well-posedness are shown with judicious choices of $\statevec u$ and $\boldsymbol\phi$. To find under what conditions a constant state is preserved, set $\mathbf u = \mathbf c = \text{constant}$ in (\ref{eq:WeakForm1}) to see that
\begin{equation}
\left( {{{\mathbf{u}}_t},\phi } \right) + \left( {\left( {\nabla  \cdot \spacevec{ \mathcal{A}} } \right){\mathbf{c}},\phi } \right) = 0,
\end{equation}
from which it follows that $\mathbf{u}_{t}=0$ if ${\nabla  \cdot \spacevec{ \mathcal{A}} }=0$.  

Global conservation is shown by selectively choosing each component of the state vector $\boldsymbol \phi$ in (\ref{eq:WeakForm2}) to be unity and again noting that the coefficient matrices are symmetric to see that the terms in the braces cancel to leave
\begin{equation}\frac{d}{dt}\int_\Omega  {{\mathbf{u}}dxdydz}  =  - \int_{\partial \Omega } {\spacevec {\mathbf{f}}  \cdot \hat ndS} .\end{equation}

To find conditions under which the initial boundary value problem is well-posed, we choose $\boldsymbol\phi = \statevec u$ in (\ref{eq:WeakForm1}) and note that
\begin{equation}
\iprod{ \nabla  \cdot  \spacevec{\mathbf{f}} + \spacevec{\mmatrix{ A}}\cdot\nabla \statevec u,\statevec u}  = \int_\Omega  {\nabla  \cdot \left( {{\statevec u^T}\spacevec{\mmatrix{ A}} \statevec u} \right)d\spacevec x}.
\end{equation}
Replacing those terms in (\ref{eq:WeakForm1}) and multiplying the equation by two gives
\begin{equation}
\frac{d}{{dt}}{\inorm{\statevec u}^2} + \int\limits_\Omega {\nabla\cdot \left( {{{\statevec{u}}^T}\spacevec{\mmatrix{A}}\statevec{u}} \right)d\spacevec{x} } + {\iprod{\left(\nabla\cdot\spacevec{\mmatrix{A}}\right)\statevec u,\statevec u}}= 0.
\end{equation}
Gauss' theorem allows us to replace
the second term by a surface integral so
\begin{equation}
\frac{d}{{dt }}\inorm{\statevec u}^2 + \int\limits_{\partial \Omega } {{{\statevec{u}}^T}\spacevec{\mmatrix{A}} \cdot \hat n{\statevec u}dS}=- {\iprod{{\left(\nabla\cdot\spacevec{\mmatrix{A}}\right){\statevec u}},{\statevec u}}} .
\end{equation}
We bound the right hand side by
\begin{equation}
-{\iprod{{\left(\nabla\cdot\spacevec{\mmatrix{A}}\right){\statevec u}},\statevec u}}  \leqslant \mathop {\max }\limits_\Omega  \left|\left| \nabla\cdot\spacevec{\mmatrix{A}} \right|\right|_{2}{\inorm{\mathbf u}^2} \equiv 2\gamma {\inorm{\statevec u}^2}
 \end{equation}
 so
 \begin{equation}
\frac{d}{{dt }}\left(e^{-2\gamma t}\inorm{\statevec u}^2\right)\le e^{-2\gamma t}\int\limits_{\partial \Omega } {{{\statevec{u}}^T}\spacevec{\mmatrix{A}} \cdot \hat n{\statevec u}dS}.
\end{equation}
Integrating over the time interval $[0,T]$ we write the energy in terms of the initial value and a boundary integral
\begin{equation}
\inorm{\statevec{u}(T)}^{2}\le e^{2\gamma T}\inorm{\statevec{u}(0)}^{2}+\int_{0}^{T} {\int_{\partial \Omega } {e^{2\gamma(T-t)}{{\statevec{u}}^T}\spacevec{\mmatrix{A}} \cdot \hat n{\statevec u}dS}dt}.
\end{equation}

To properly pose the problem we must impose appropriate boundary conditions. From (\ref{eq:mDHyperbolicConditions_DAK}), we separate the waves traveling to the left and right of the boundary relative
to $\hat n$ as
\begin{equation}
\spacevec{\mmatrix{A}} \cdot \hat n = \sum\limits_{m = 1}^3 {{\mmatrix A^{(m)}}{{\hat n}_m}}  = \mmatrix R\Lambda {\mmatrix R^{ - 1}} = \mmatrix P{\Lambda ^ + }{\mmatrix R^{ - 1}} + \mmatrix R{\Lambda ^ - }{\mmatrix R^{ - 1}} \equiv {\mmatrix A^ + } + {\mmatrix A^ - },
\end{equation}
where $\Lambda^{\pm} = \Lambda \pm \left|\Lambda\right|$ and we have left off the explicit dependence on $\hat n$.
When we replace the values of $\statevec{u}$ along the boundary associated with the incoming $\Lambda^{-}$ waves with a boundary state, $\statevec{g}$, the solution can be bounded in terms of the initial and boundary data,
\begin{equation}
\begin{split}
\inorm{\statevec{u}(T)}^{2} +\int_{0}^{T} {\int_{\partial \Omega } {{{\statevec{u}}^T}{\mmatrix{A}^{+}}{\statevec u}dS}dt}
&\le e^{2\gamma T}\inorm{\statevec{u}(0)}^{2}+\int_{0}^{T} {\int_{\partial \Omega } {e^{2\gamma(T-t)}{{\statevec{g}}^T}\left|\mmatrix{A}^{-}\right|{\statevec g}dS}dt} 
\\&\le e^{2\gamma T}\left\{\inorm{\statevec{u}(0)}^{2}+\int_{0}^{T} {\int_{\partial \Omega } {{{\statevec{g}}^T}\left|\mmatrix{A}^{-}\right|{\statevec g}dS}dt}\right\}.
\end{split}
\label{eq:dDimensionWellPosednessLinear_DAK}
\end{equation}
Furthermore, if $\nabla\cdot\spacevec{\mathcal{A}} = 0$, $\gamma = 0$ and the energy does not grow in time except for energy introduced at the boundaries,
\begin{equation}
\inorm{\statevec{u}(T)}^{2} +\int_{0}^{T} {\int_{\partial \Omega } {{{\statevec{u}}^T}{\mmatrix{A}^{+}}{\statevec u}dS}dt}
\le \inorm{\statevec{u}(0)}^{2}+\int_{0}^{T} {\int_{\partial \Omega } {{{\statevec{g}}^T}\left|\mmatrix{A}^{-}\right|{\statevec g}dS}dt} .
\label{eq:dDimensionWellPosednessLinearG=0}
\end{equation}

\section{A Polynomial Spectral Calculus}
\label{sec:PolyCalc}
To follow the continuous problem analysis as closely as possible, we introduce a discrete calculus that looks and behaves 
like the continuous one as much as possible. We define the calculus for the reference domain $E=[-1,1]^{3}$ with coordinates $\spacevec \xi = \left(\xi,\eta,\zeta\right) = \xi \hat \xi + \eta\hat\eta + \zeta\hat\zeta = \sum\limits_{m = 1}^3 {{\xi ^{(m)}}{{\hat \xi }^m}}$. Corresponding forms hold for two dimensional problems.

We represent functions of the reference domain coordinates by polynomials of degree $N$ or less, i.e. as elements of $\mathbb{P}^{N}(E)\subset\mathbb{L}^{2}(E)$. A basis for the polynomials on $E$ is the tensor product of the one dimensional Lagrange basis. 
Using that basis, we write a polynomial, $U$, in terms of nodal values ${U}_{ijk}={U}\left(\xi_{i},\eta_{j},\zeta_{k}\right)$ as an upper case letter, which for three space dimensions is
\begin{equation}{{U}} = \sum\limits_{i,j,k = 0}^N {{{{U}}_{ijk}}{\ell _i}(\xi){\ell _j}(\eta){\ell _k}(\zeta)}, \end{equation}
where 
\begin{equation}{\ell _l}\left( s  \right) = \prod\limits_{i = 0;i \ne l}^N {\frac{{s  - {s _i}}}{{{s _l} - {s _i}}}} \end{equation}
is the one-dimensional Lagrange interpolating polynomial with the property $\ell_{l}\left(s_{m}\right) = \delta_{lm}$, $l,m=0,1,2,\ldots,N$. The points $s_{i},\; i=0,1,2,\ldots,N$ are the interpolation points, whose locations are chosen below. We also write the interpolation operator, $\mathbb{I}^{N}:\mathbb{L}^{2}\rightarrow\mathbb{P}^{N}$, which projects square integrable functions on $E$ onto polynomials, as
\begin{equation}{\mathbb{I}^N}\left(  u \right) = \sum\limits_{i,j,k = 0}^N {{ u_{ijk}}{\ell _i}(\xi){\ell _j}(\eta){\ell _k}(\zeta)}. \end{equation}
The use of the tensor product means that one and two dimensions are special cases of three dimensions, which is why we concentrate on three dimensional geometries here.

 Derivatives of polynomials on $E$ evaluated at the nodes can be represented by matrix-vector multiplication. For instance,
 \begin{equation}
 {\left. {\frac{{\partial U}}{{\partial \xi}}} \right|_{nml}} = \sum\limits_{i,j,k = 0}^N {{{{U}}_{ijk}}{\ell' _i}({\xi_n}){\ell _j}({\eta_m}){\ell _k}({\zeta_l})}  = \sum\limits_{i = 0}^N {{U_{iml}}{{\ell '}_i}({\xi_n})}  \equiv \sum\limits_{i = 0}^N {{U_{iml}}{\mathcal{D}_{ni}}},
  \end{equation}
  where $\mathcal{D}$ is the derivative matrix. The gradient and divergence of a polynomial in three space dimensions evaluated at a point $\left( \xi_{n},\eta_{m},\zeta_{l}\right)$ are therefore
  \begin{equation}\begin{gathered}
  {\left. {\nabla  U} \right|_{nml}} = \sum\limits_{i = 0}^N {{U_{iml}}{\mathcal{D}_{ni}}} \hat \xi + \sum\limits_{j = 0}^N {{U_{njl}}{\mathcal{D}_{mj}}} \hat \eta + \sum\limits_{k = 0}^N {{U_{nmk}}{\mathcal{D}_{lk}}} \hat \zeta, \hfill \\ 
  {\left. {\nabla \cdot\spacevec F} \right|_{nml}} = \sum\limits_{i = 0}^N {{F^{(1)}_{iml}}{\mathcal{D}_{ni}}}  + \sum\limits_{j = 0}^N {{F^{(2)}_{njl}}{\mathcal{D}_{mj}}}  + \sum\limits_{k = 0}^N {{F^{(3)}_{nmk}}{\mathcal{D}_{lk}}}.  \hfill \\
\end{gathered} 
  \end{equation}

The use of the calculus that we develop depends on
the choice that the interpolation nodes, $s_{i}$, are the nodes of the Legendre-Gauss-Lobatto (LGL) quadrature.
We represent the one dimensional LGL quadrature of a function $g(s)$ as
\begin{equation}
\int_{ - 1}^1 {gds}\approx \sum\limits_{i = 0}^N {{g\left(s_{i}\right)}{\omega _i}}\equiv \int_N {gds}    ,
 \end{equation} 
 where the $\omega_{i}$ are the LGL quadrature weights.
 The quadrature is exact if $g\in\mathbb{P}^{2N-1}$. By tensor product extension, we write three dimensional volume integral approximations as
 \begin{equation}
 \int_{E,N} {gd \xi d\eta d\zeta}  \equiv \sum\limits_{i,j,k = 0}^N {{g_{ijk}}{\omega_{ijk}}},
  \end{equation}
 where $\omega_{ijk}=\omega_{i}\omega_{j}\omega_{k}$. Two-dimensional surface integral approximations are
 \begin{equation}
 \begin{split}
 \int_{\partial E,N} {\spacevec g \cdot \hat ndS}  &= \sum\limits_{i,j = 0}^N {\left. {{\omega_{ij}}{g^{(1)}}\left( {\xi,{\eta_i},{\zeta_j}} \right)} \right|_{\xi =  - 1}^1}  + \sum\limits_{i,j = 0}^N {\left. {{\omega_{ij}}{g^{(2)}}\left( {{\xi_i},\eta,{\zeta_j}} \right)} \right|_{\eta =  - 1}^1}  + \sum\limits_{i,j = 0}^N {\left. {{\omega_{ij}}{g^{(3)}}\left( {{\xi_i},{\eta_j},\zeta} \right)} \right|_{\zeta =  - 1}^1} 
\\&\equiv\int_N {\left. {{g^{(1)}}d\eta d\zeta } \right|} _{\xi  =  - 1}^1 + \int_N {\left. {{g^{(2)}}d\xi d\zeta } \right|} _{\eta  =  - 1}^1 + \int_N {\left. {{g^{(3)}}d\xi d\eta } \right|} _{\zeta  =  - 1}^1.
 \end{split}
 \end{equation}
 Two space dimensional areas and edge integrals are defined similarly. 
 
We define the discrete inner product of two functions $f$ and $g$ and the discrete norm of $f$ from the quadrature
 \begin{equation}{\left( {f,g} \right)_{E,N}} = \int_{E,N} {fgd \xi d \eta d \zeta }  \equiv \sum\limits_{i,j,k = 0}^N {{f_{ijk}}{g_{ijk}}{\omega_{ijk}}} ,\quad {\left\| f \right\|_{E,N}} = \sqrt {{{\left( {f,f} \right)}_{E,N}}}. \end{equation}
 The definition is extended for vector arguments like
 \begin{equation}\spacevec{\statevec f} = \sum\limits_{m = 1}^3 {{\statevec f^{(m)}}{{\hat \xi }^m}}, \end{equation}
 for a state vector $\statevec f^{(m)} = [f^{(m)}_{1}\;f^{(m)}_{2}\;\ldots \; f^{(m)}_{p}]^{T}$
  as
 \begin{equation}{\left( {\spacevec{\statevec f},\spacevec {\statevec g}} \right)_N} = \int_{E,N} {\sum\limits_{m = 1}^3 {{{\left( {{\statevec f^{(m)}}} \right)}^T}{\statevec g^{(m)}}d \xi d \eta d \zeta } }  = \sum\limits_{i,j,k = 0}^N {{\omega_{ijk}}\sum\limits_{m = 1}^3 {{{\left( {\statevec f_{ijk}^{(m)}} \right)}^T}\statevec g_{ijk}^{(m)}} }, \end{equation}
 and similarly for other arguments.
 
 The Lagrange basis functions are orthogonal with respect to the discrete inner product \cite{CHQZ:2006}. In one space dimension, for instance, $\left(\ell_{i},\ell_{j}\right)_{E,N} = \omega_{j}\delta_{ij}$.
 Also, from the definitions of the interpolation operator and the discrete inner product,
 \begin{equation}
 {\left( {f,g} \right)_{E,N}} = {\left( {{\mathbb{I}^N}\left( f \right),{\mathbb{I}^N}\left( g \right)} \right)_{E,N}}.
 \end{equation} 
Finally, the discrete norm is equivalent to the continuous norm \cite{ISI:A1982NE30900005} in that for $U\in\mathbb{P}^{N}$,
 \begin{equation}
 {\left\| U \right\|_E} \leqslant {\left\| U \right\|_{E,N}} \leqslant C{\left\| U \right\|_E},
 \end{equation}
 where $C$ is a constant.

The crucial property for the analysis of the discrete approximation is the \emph{summation by parts} (SBP) property satisfied by the LGL quadrature. Let $U,V\in\mathbb{P}^{N}$. Then exactness of the LGL quadrature implies that
\begin{equation}
{\int_N {UV'dx}  = \left. {UV} \right|_{ - 1}^1 - \int_N {U'Vdx} }\quad(Summation\; By\; Parts).
\label{eq:SumByParts_DAK}
\end{equation}
The result extends to all space dimensions \cite{gassner2010} with
\begin{equation}
\begin{gathered}
  {\left( {{U_\xi },V} \right)_N} = \int_N {\left. {UVd\eta d\zeta } \right|} _{\xi  =  - 1}^1 - {\left( {U,{V_\xi }} \right)_N} \hfill \\
  {\left( {{U_\eta },V} \right)_N} = \int_N {\left. {UVd\xi d\zeta } \right|} _{\eta  =  - 1}^1 - {\left( {U,{V_\eta }} \right)_N} \hfill \\
  {\left( {{U_\zeta },V} \right)_N} = \int_N {\left. {UVd\xi d\eta } \right|} _{ \zeta=  - 1}^1 - {\left( {U,{V_\zeta }} \right)_N}. \hfill \\ 
\end{gathered} 
\label{eq:SBPDerivative}\end{equation}

We can use (\ref{eq:SumByParts_DAK}) and (\ref{eq:SBPDerivative}) to formulate a discrete integral calculus. If we replace $U$ in (\ref{eq:SBPDerivative}) by the components of a vector $\spacevec F$, and sum, we get the \emph{Discrete Extended Gauss Law (DXGL)} originally derived in \cite{gassner2010}: For any vector of polynomials $\spacevec F\in\mathbb{P}^{N}$ and any polynomial $V\in\mathbb{P}^{N}$,
 \begin{equation}
 {\iprodN{\nabla  \cdot \spacevec F,V} = \int_{\partial E ,N} {\spacevec F \cdot \hat nVdS}  - \iprodN{\spacevec F,\nabla V}}\quad (Discrete \; Extended \; Gauss\; Law),
 \label{eq:DiscreteGreens_DAK}
 \end{equation}
  where $\hat n$ is the unit outward normal at the faces of $E$.
 Carrying this further, if we set $V=1$ we get the \emph{Discrete Gauss Law} (DGL)
  \begin{equation}
 \iprodN{\nabla  \cdot \spacevec F,1} =\int_{ E ,N} {\nabla\cdot\spacevec Fd \xi d\eta d\zeta} = \int_{\partial E ,N} {\spacevec F \cdot \hat ndS} \quad(Discrete\; Gauss\; Law).
 \label{eq:DiscreteDivergence_DAK}
 \end{equation}
The DGL is exact for polynomial arguments. By using the appropriate definitions for the inner products, both discrete Gauss laws extend to hold for state vectors $\spacevec {\statevec F}$ and $\statevec V$. 

Next, we see that if we replace the vector flux $\spacevec F$ in (\ref{eq:DiscreteGreens_DAK}) with $\nabla \Phi \in\mathbb{P}^{N}$, then we get the discrete version of Green's first identity,
\begin{equation}
{\left( {{\nabla ^2}\Phi ,V} \right)_N} + {\left( {\nabla \Phi ,\nabla V} \right)_N} = \int_{\partial E,N} {\nabla \Phi  \cdot \hat nVdS} \quad(Discrete\; Green's\;First\; Identity).
\end{equation}
Swapping the variables $\Phi$ and $V$ and subtracting from the original gives Green's second identity
\begin{equation}
{\left( {{\nabla ^2}\Phi ,V} \right)_N} - {\left( {{\nabla ^2}V,\Phi } \right)_N} = \int_{\partial E,N} {\left( {\nabla \Phi  \cdot \hat nV - \nabla V \cdot \hat n\Phi } \right)dS} \quad(Discrete\; Green's\;Second\; Identity)
 \end{equation}
The discrete Green's identities would be useful to prove stability of continuous Galerkin spectral element methods of second order problems.
 
 Other identities that do not involve quadratic products of polynomial arguments hold discretely through exactness of the LGL quadrature. For instance, 
 \begin{equation}
 \int_{E,N} {\nabla Vd \xi d\eta d\zeta }  = \int_{\partial E,N} {V\hat n dS} 
 \end{equation}
 and
 \begin{equation}
 \int_{E,N} {\nabla  \times \spacevec Fd \xi d\eta d\zeta }  = \int_{\partial E,N} {\hat n \times \spacevec FdS}. 
 \end{equation}
 What we see, then, is that the well-known integral identities hold due to either integration or summation by parts. 

 Whereas integration rules hold discretely, product differentiation rules do not usually hold because differentiation and interpolation do not always commute. For instance, the product rule does not generally hold. That is, for polynomials $U,V$,
 \begin{equation}
 \nabla \left({\mathbb{I}^N}\left( {UV}\right) \right) \ne {\mathbb{I}^N}\left( {U\nabla V} \right) + {\mathbb{I}^N}\left( {V\nabla U} \right)
 \label{eq:ProdRule}
 \end{equation}
 unless the product $UV\in\mathbb{P}^{N}$.
 Differentiation and interpolation do not commute because of aliasing errors that arise from projecting the product onto a polynomial of degree less than or equal to $N$ \cite{CHQZ:2006}. 
 
\section{Discontinuous Galerkin Spectral Element Approximations}
We now use the polynomial calculus introduced in (\ref{sec:PolyCalc}) to formulate and analyze discontinuous Galerkin spectral element approximations 
in three space dimensions. The steps to derive two dimensional approximations are identical. 
The domain $\Omega$ is subdivided into $N_{el}$ nonoverlapping hexahedral elements, $e^{r}, r=1,2,\ldots,N_{el}$. We assume here that the subdivision is conforming. Each element is mapped from the reference element $E$ by a transformation $\spacevec x = \spacevec X\left(\spacevec \xi\right)$. 
From the transformation, we define the three covariant basis vectors
\begin{equation}
\spacevec{a}_{i}=\partialderivative{\spacevec{X}}{\xi^{i}}\quad i=1,2,3,
\end{equation}
and (volume weighted) contravariant vectors, formally written as
\begin{equation}
\mathcal{J}\spacevec{a}^{i}=\spacevec{a}_{j}\times\spacevec{a}_{j}, \quad (i,j,k)\;\text{cyclic},
\end{equation}
where
\begin{equation}
\mathcal{J}=\spacevec{a}_{i}\cdot\left(\spacevec{a}_{j}\times \spacevec{a}_{k}\right),\quad (i,j,k)\;\text{cyclic}
\end{equation}
is the Jacobian of the transformation. 

Under the mapping, the 
divergence of a spatial vector flux can be written compactly in terms of the reference space variables
as
\begin{equation}
\nabla \cdot \spacevec {\statevec f} 
= \frac{1}{\mathcal{J}}\sum\limits_{i = 1}^3 {\frac{\partial }{{\partial {\xi ^i}}}\left( {\mathcal{J}{{\spacevec a}^i} \cdot \spacevec{ \statevec f}} \right)}
= \frac{1}{\mathcal{J}}\sum\limits_{i = 1}^3 {\frac{\partial \tilde {\statevec f}^{i}}{{\partial {\xi ^i}}}}
= \frac{1}{\mathcal{J}}{\nabla _\xi } \cdot \contravec{ {\statevec f}}.
\label{eq:CompSpaceDivergence_DAK}
\end{equation}
The vector $\contravec{{\statevec f}}$ is the volume weighted contravariant flux whose components are $\contravec{{\statevec f}}^{i}={\mathcal{J}{{\spacevec a}^i} \cdot \spacevec{ {\statevec f}}}$.

The conservation law is then represented on the reference domain by another conservation law
\begin{equation}
{\mathcal{J}\statevec u_t} +\nabla_{\xi}\cdot\left(\contravec{\mmatrix{A}}\statevec u\right)  = 0,
\end{equation}
where we have defined the (volume weighted) contravariant coefficient matrices
\begin{equation}
{\mmatrix A^i} = \mathcal{J}{{\spacevec a}^i} \cdot \spacevec{\mmatrix{ A}}
\end{equation}
and
\begin{equation}\contravec{\mmatrix{ A}} = \sum\limits_{i = 1}^3 {{\mathcal A^i}{{\hat \xi }^i}}. \end{equation}
We can also construct the nonconservative form of the system on the reference domain using the chain rule,
\begin{equation}
\mathcal{J}{\statevec u_t} + \left(\nabla_{\xi}\cdot\contravec{\mmatrix{A}}\right)\statevec u  +\contravec{\mmatrix{A}}\cdot\nabla_{\xi}\statevec u  = 0.
\end{equation}

We construct weak forms of the conservative and nonconservative equations by taking the inner product of the equations with a test function $\boldsymbol\phi \in \mathbb{L}^{2}(E)$ and applying extended Gauss Law to the space derivative terms,
\begin{equation}
\iprod{\mathcal{J}\statevec{u}_{t},\boldsymbol\phi}_{E}+\int_E {\contravec{\mathbf{ f}} \cdot \hat n^{T}\boldsymbol\phi dS}-\iprod{\contravec{\mathbf{f}},\nabla_{\xi}\boldsymbol\phi}_{E}=0
\label{eq:WeakConservative}
\end{equation}
and
\begin{equation}
\iprod{\mathcal{J}{\statevec u_t},\boldsymbol\phi} +\int_E {\contravec{\mathbf{ f}} \cdot \hat n^{T}\boldsymbol\phi dS}- \iprod{\statevec u,\nabla_{\xi}\cdot\contravec{\statevec{f}}\left(\boldsymbol\phi\right)}_{E}+ \iprod{\left(\nabla_{\xi}\cdot\contravec{\mmatrix{A}}\right)\statevec u,\boldsymbol\phi}_{E}    = 0.
\label{eq:WeakNonConservative}
\end{equation}
When we average the two equations (\ref{eq:WeakConservative})
and (\ref{eq:WeakNonConservative}) we get the split weak form
\begin{equation}
\iprod{\mathcal{J}{\statevec u_t},\boldsymbol\phi}    - \frac{1}{2}\left\{ {\iprod{\contravec{\statevec{f}}(\statevec u),\nabla_{\xi} \boldsymbol \phi}_{E} + \iprod{\statevec u,\nabla_{\xi}  \cdot \contravec{\statevec{f}}\left( \boldsymbol\phi  \right)}_{E} - \iprod{\left(\nabla_{\xi}  \cdot \contravec{\mmatrix A}\right)\statevec u,\boldsymbol\phi }}_{E} \right\} + \int_{\partial E} {\left(\contravec{\statevec{f}} \cdot {{\hat n} }\right)^{T}\boldsymbol\phi d S}=0.
\label{eq:GeneralWeakFormLinear_DAK}
\end{equation}
\subsection{The DGSEM}
The original DG spectral element method introduced by Black \cite{Black1999} starts from the conservative weak form (\ref{eq:WeakConservative}). We use the calculus now to show that it is stable if the coefficient matrices $\contravec{\mathcal A}$ are constant. If, in addition, characteristic boundary conditions are used at physical boundaries, the approximation is optimally stable in the sense that the global energy discretely matches (\ref{eq:dDimensionWellPosednessLinearG=0}).

To construct the approximation,
one approximates the solutions, fluxes, coefficient matrices and Jacobian with polynomial interpolants on element $e^{r}\rightarrow E$ by
\begin{equation}\begin{gathered}
  {\mathbf{u}} \approx {\mathbf{U}}^{r} \in {\mathbb{P}^N} \hfill \\
  \contravec{\mathbf{f}} \approx \contravec{\mathbf{F}}^{r}\left(\statevec U\right) = {\mathbb{I}^N}\left( {{\mathbb{I}^N}\left( {\tilde{ \mathcal{A}}} ^{}\right){\mathbf{U}}} \right) = \sum\limits_{i,j,k = 0}^N {{{\tilde{ \mathcal{A}}}_{ijk}}{{\mathbf{U}}_{ijk}}{\ell _i}\left( \xi  \right){\ell _j}\left( \eta  \right){\ell _k}\left( \zeta  \right)}  \hfill \\
  {\tilde{ \mathcal{A}}} \approx {\mathbb{I}^N}\left( {\tilde{ \mathcal{A}}} \right) \hfill \\
  \mathcal{J}^{r} \approx J^{r} = {\mathbb{I}^N}\left( \mathcal{J} ^{r}\right). \hfill \\ 
\end{gathered} 
\label{eq:Approxes}
\end{equation}
From this point, we leave off the superscripts $r$ and subscripts $\xi$ on $\nabla_{\xi}$ unless necessary.

To continue the contruction, one replaces the continuous inner products by the discrete inner products, here being Gauss-Lobatto quadratures. The normal boundary flux is replaced by a consistent numerical flux, $\contravec{\statevec{f}}\leftarrow \contravec{\statevec{F}}^{*}\left(\statevec U^{L},\statevec U^{R};\hat n\right)$ where $\statevec U^{L,R}$ are the left and right states at the element boundary, measured with respect to the outward normal, $\hat n$. The numerical flux ensures continuity of the normal flux at element faces. Finally, $\boldsymbol \phi$ is restricted to elements of $\mathbb{P}^{N}$.  
The result of the approximations is the formal statement of the method
\begin{equation}
 [DGSEM]\quad {\left( {J{\mathbf{U}}_{t},\boldsymbol\phi } \right)_N}+ \int_{\partial E,N} {{{\tilde {\mathbf{F}} }^{*,T}}\boldsymbol\phi dS}   - {{\left( {{\contravec {\mathbf{F}} }\left( {\mathbf{U}} \right),{\nabla  }\boldsymbol\phi } \right)}_N} = 0.
 \label{eq:DGSEM-W}
\end{equation}
Details for going from the formal statement to the form to implement can be found in \cite{Kopriva:2009nx}. 

Alternate, yet algebraically equivalent forms of the DGSEM can be derived by applying the DXGL. For instance, if we apply the DXGL to the last inner product in (\ref{eq:DGSEM-W}) we get the algebraically equivalent form
\begin{equation}
{\left( {J{\mathbf{U}}_{t},\boldsymbol\phi } \right)_N} + \int_{\partial E,N} {{{\left\{ {{{\tilde {\mathbf{F}} }^*} - {\tilde {\mathbf{F}} } \cdot \hat n} \right\}}^T}\boldsymbol\phi dS}  + {{\left( {{\nabla  } \cdot \tilde {\statevec{F}}\left( {\mathbf{U}} \right),\boldsymbol\phi } \right)}_N}= 0.
 \label{eq:DGSEM-S}
\end{equation}

If the contravariant coefficient matrices are constant, implying that the original problem is constant coefficient and the elements are rectangular in shape, then the DGSEM approximation is strongly stable. To show stability, we set $\boldsymbol\phi = \statevec{U}$ in (\ref{eq:DGSEM-S})
and define the volume weighted norm
\begin{equation}
\left\| {\mathbf{U}} \right\|_{J,N}^2 \equiv {\left( {J{\mathbf{U}},{\mathbf{U}}} \right)_N}.
\end{equation}
Then
\begin{equation}
\frac{1}{2}\frac{d}{{dt}}\left\| \statevec U \right\|_{J,N}^2 + \int_{\partial E,N} {{{\left\{ {{{\tilde {\mathbf{F}} }^*} - {\tilde {\mathbf{F}} } \cdot \hat n} \right\}}^{T}}\statevec U dS} + {{\left( {{\nabla  } \cdot \tilde {\statevec{F}}\left( {\mathbf{U}} \right),\statevec U } \right)}_N} = 0.
\label{eq:DGSEMStab1}
\end{equation}
With constant coefficients, the volume term in (\ref{eq:DGSEMStab1}) can be converted to a surface quadrature. The coefficient matrices being constant and symmetric implies that
\begin{equation}
\begin{split}
{\left( {\nabla  \cdot \contravec{\statevec F}\left(\statevec  U \right),\statevec U} \right)_N} &= \underbrace{{\left( {\nabla  \cdot \mathbb{I}^{N}\left( {\contravec{\mathcal{A}}\statevec U} \right),\statevec U} \right)_N} = {\left( {\contravec{\mathcal{A}} \cdot \nabla \statevec U,\statevec U} \right)_N}}_{*} = {\left( {\nabla \statevec U,\contravec{\mathcal{A}}\statevec U} \right)_N} \\&= {\left( {\nabla \statevec U,\contravec{\statevec F}\left( \statevec U \right)} \right)_N}.
\end{split}
\label{eq:KeyStep}
\end{equation}
The key step is the second marked with the ``*'', where the product rule applies because $\contravec{\mathcal A}\statevec U\in\mathbb{P}^{N}$ when $\contravec{\mathcal A}$ is constant.
We then substitute the equivalence into the DXGL to see that
\begin{equation}
{\left( {\nabla  \cdot \contravec{\statevec F}\left( \statevec U \right),\statevec U} \right)_N} = \frac{1}{2}\int_{\partial E,N} {{{\left( {\contravec{\statevec F} \cdot \hat n} \right)}^T}\statevec UdS} .
\end{equation}
Therefore, the local energy changes according to
\begin{equation}
\frac{1}{2}\frac{d}{{dt}}\left\| \statevec U \right\|_{J,N}^2 + \int_{\partial E,N} {{{\left\{ {{{\tilde {\mathbf{F}} }^*} - \frac{1}{2}{\tilde {\mathbf{F}} } \cdot \hat n} \right\}}^{T}}\statevec U dS}  = 0,
\label{eq:DGSEMLocalEnergyChange}
\end{equation}
and stability depends solely on what happens on the element faces.

The change in the total energy is found by summing over all the elements. Although the numerical flux is continuous at element interfaces, the solution and flux are discontinuous. If we define the jump in a quantity with the usual notation $\jump{V}=V^{R}-V^{L}$, then
\begin{equation}
\begin{split}
\frac{d}{{dt}}\left( {\sum\limits_{r = 1}^{N_{el}} {\inorm{\statevec U^{r}}_{J,N}^2} } \right) 
\leqslant 
-2
&\left\{ \sum\limits_{Boundary\atop Faces} {\int_{\partial E,N} {{{\left( {{{\statevec F}^*} - \frac{1}{2}\statevec F \cdot \hat n} \right)}^T}\statevec UdS} }  \right.
\\&- \left. \sum\limits_{Interior\atop Faces} {\int_{\partial E,N} {\left( {{{\statevec F}^{*,T}}\jump{\statevec U} - \frac{1}{2}\jump{{{\left( {\statevec F \cdot \hat n} \right)}^T}\statevec U}} \right)dS} }  \right\}.
\end{split}
\label{eq:IntFaceTerms1}
\end{equation}

Stability is determined, therefore, only by the influence of the jumps at the element boundaries and the physical boundary approximations through the numerical flux.
 For linear problems, it is natural to choose an upwinded or central flux,
 \begin{equation}
 {\contravec{\statevec F}^*}\left( {{{\mathbf{U}}^L},{{\mathbf{U}}^R};\hat n} \right) = \frac{1}{2}\left\{ {\contravec{\statevec F}\left( {{{\mathbf{U}}^L}} \right) \cdot \hat n + \contravec{\statevec F}\left( {{{\mathbf{U}}^R}} \right) \cdot \hat n} \right\} - \sigma \frac{{\left| {\tilde {\mathcal A} \cdot \hat n} \right|}}{2}\left\{ {{{\mathbf{U}}^R} - {{\mathbf{U}}^L}} \right\},
 \label{eq:RiemannFlux}
 \end{equation}
 where $\sigma=0$ is the central flux and $\sigma = 1$ is the fully upwind flux. With this flux \cite{Kopriva:2014yq},
 \begin{equation}
{{{{\contravec {\statevec F}}^{*,T}}\jump{\statevec U} - \frac{1}{2}\jump{{{\left( {\contravec {\statevec F} \cdot \hat n} \right)}^T}\statevec U}} } = -\frac{\sigma}{2}\jump{\statevec U}^{T}{{\left| {\tilde {\mathcal A} \cdot \hat n} \right|}}\jump{\statevec{U}}\le 0,
\end{equation}
so that the interior face terms in (\ref{eq:IntFaceTerms1}) are dissipative. To match the PDE energy bound, (\ref{eq:dDimensionWellPosednessLinearG=0}), the fully upwind flux must be used at the physical boundaries.
With exterior values $\statevec g$ set along incoming characteristics \cite{David-A.-Kopriva:2015dp} and when $\sigma=1$,
\begin{equation}
{\left( {{{ {\statevec{F}} }^*} - \frac{1}{2} \contravec{\statevec{F}}  \cdot \hat n} \right)^T}{\statevec{U}} = \frac{1}{2}{{\statevec{U}}^T}{\mmatrix A^ + }{\statevec{U}} + \frac{1}{2}\left\| {\sqrt {\left| {{\mmatrix A^ - }} \right|} {\statevec{U}} - \sqrt {\left| {{\mmatrix A^ - }} \right|} {{\statevec{g}} }} \right\|_2^2 - \frac{1}{2}{\statevec{g}} ^T\left| {{\mmatrix A^ - }} \right|{{\statevec{g}} }.
\label{eq:PhyEnergydissipation}
\end{equation}
If we define the total energy by
\begin{equation}
\inorm{\statevec U}_{J,N}^{2} = \sum\limits_{r = 1}^K {\inorm{\statevec U^{r}}_{J,N}^2},
\end{equation}
and integrate (\ref{eq:IntFaceTerms1})in time, the total energy satisfies (c.f. (\ref{eq:dDimensionWellPosednessLinear_DAK}))
\begin{equation}
\begin{split}
\inorm{\statevec U(T)}_{J,N}^2 + \sum\limits_{Boundary\atop Faces} {\int_{0}^{T}\int_{\partial E,N} {{{\statevec U}^T}{\mmatrix A^ + }\statevec UdSdt} }  
\leqslant \inorm{\statevec U(0)}_{J,N}^2 + \sum\limits_{Boundary\atop Faces} {\int_{0}^{T}\int_{\partial E,N} {\statevec g ^T\left| {{\mmatrix A^ - }} \right|{{\statevec g} }dSdt} } .
\end{split}
\label{eq:MDStabilityLinearNoGrowth_DAK}
\end{equation}
Finally, if the interpolant of the Jacobian is bounded from below, $J>0$, then for some positive constants $c$ and $C$ \cite{Kopriva:2014yq},
\begin{equation}
c\left\| {{\vec U}} \right\|_{{L^2}\left( \Omega \right)}^2 \leqslant \left\| {{\vec U}} \right\|_{J,N}^2 \leqslant C\left\| {{\vec U}} \right\|_{{L^2}\left( \Omega \right)}^2,
\label{eq:NormEquivalence}
\end{equation}
which says that, like the continuous solution, the energy approximate solution is bounded by the data in the continuous norm over the entire domain.

\subsection{Stabilization by Split Form}
If the contravariant coefficient matrices are not constant, then the key step in (\ref{eq:KeyStep}) does not hold because interpolation and differentiation do not commute. We show now that stability hangs on whether or not the dissipation introduced by the numerical flux at the element interfaces and by the characteristic boundary conditions is sufficient to counterbalance the aliasing errors associated with the volume term that remains. That balance shows why the approximation [DGSEM] can be, but does not have to be, stable for variable coefficient problems or curved elements.

The use of the polynomial calculus allows us to quickly and compactly construct four algebraically equivalent representations of a split form approximation \cite{Kopriva:2014yq} that is strongly stable, constant state preserving and globally conservative for non-constant coefficient problems where the coefficient variation is due to inherent variability in the PDE and/or due to variability introduced by the coefficient mappings from curved elements to the reference element. It also allows us to simplify the analysis done, for example, in \cite{Kopriva2016274}.

The result of applying the approximations (\ref{eq:Approxes}) and LGL quadrature to (\ref{eq:GeneralWeakFormLinear_DAK}) is the first split form of the DGSEM used in \cite{Kopriva:2014yq}. In accordance to common terminology, this is the ``weak'' form
\begin{equation}
\begin{split}
\left[ W \right]\quad {\left( {J{\mathbf{U}}_{t},\boldsymbol\phi } \right)_N}   &- \frac{1}{2}\left\{ {{{\left( {{\contravec {\mathbf{F}} }\left( {\mathbf{U}} \right),{\nabla  }\boldsymbol\phi } \right)}_N} + {{\left( {{\mathbf{U}},{\nabla  } \cdot{\contravec {\mathbf{F}} }\left( \boldsymbol\phi  \right)} \right)}_N} - {{\left( {{\nabla  } \cdot \left({\mathbb{I}^N}\left( {\tilde {\mmatrix A}} \right)\right){\mathbf{U}},\boldsymbol\phi } \right)}_N}} \right\} \\&+ \int_{\partial E,N} {{{\tilde {\mathbf{F}} }^{*,T}}\boldsymbol\phi dS}= 0.
\end{split}
\label{eq:[W]}
\end{equation}
We get alternative, yet algebraically equivalent forms by applying the DXGL (\ref{eq:DiscreteGreens_DAK}) to selected terms in (\ref{eq:[W]}). When we apply the DXGL  to the first two inner products in the braces and use the fact that the coefficient matrices are symmetric we get the  ``strong'' form
\begin{equation}
\begin{split}
\left[ S \right]\quad {\left( {J{\mathbf{U}}_{t},\boldsymbol\phi } \right)_N}   &+ \frac{1}{2}\left\{ {{{\left( {{\nabla  } \cdot \tilde {\statevec{F}}\left( {\mathbf{U}} \right),\boldsymbol\phi } \right)}_N} + {{\left( {{\mathbb{I}^N}\left( {\tilde {\mmatrix A}} \right) \cdot {\nabla  }{\mathbf{U}},\boldsymbol\phi } \right)}_N} + {{\left( {{\nabla  } \cdot \left({\mathbb{I}^N}\left( {\tilde {\mmatrix A}} \right)\right){\mathbf{U}},\boldsymbol\phi } \right)}_N}} \right\} \\&+ \int_{\partial E,N} {{{\left\{ {{{\tilde {\mathbf{F}} }^*} - {\tilde {\mathbf{F}} } \cdot \hat n} \right\}}^T}\boldsymbol\phi dS}= 0.
\end{split}
\end{equation}
If we rearrange the terms in $[S]$ to ``strong+correction'' form
\begin{equation}
\begin{split}
\left[ {SC} \right]\quad {\left( {J{{\mathbf{U}}_t},\boldsymbol\phi } \right)_N} 
&+ 
\int_{\partial E,N} {{{\left\{ {{{\tilde {\mathbf{F}}}^*} - \tilde {\mathbf{F}} \cdot \hat n} \right\}}^T}\boldsymbol\phi dS} 
 + 
{\left( {\nabla  \cdot \contravec{\statevec F}\left( {\mathbf{U}} \right),\boldsymbol\phi } \right)_N} \\&+ \frac{1}{2}{\left( {\left\{{\mathbb{I}^N}\left( {\tilde {\mmatrix A}} \right) \cdot \nabla {\mathbf{U}} 
+
 \nabla  \cdot \left({\mathbb{I}^N}\left( {\tilde {\mmatrix A}} \right)\right){\mathbf{U}}
  - 
  \nabla  \cdot \contravec{\statevec F}\left( {\mathbf{U}} \right)\right\},\boldsymbol\phi } \right)_N} = 0,
\end{split}
\end{equation}
we see that the split form approximation is the strong form of the original DGSEM (\ref{eq:DGSEM-S}) plus a correction term in the amount by which the product rule (\ref{eq:ProdRule}) does not hold. When the product rule does hold, such as when the contravariant coefficient matrices are constant, the correction term vanishes and we are back to the original scheme of Black, [DGSEM].

We get a fourth algebraically equivalent ``directly stable'' form by applying the DXGL to only the first inner product in the braces of the weak form $[W]$,
\begin{equation}
\begin{split}
\left[ {DS} \right]\quad {\left( {J{\mathbf{U}}_{t},\boldsymbol\phi } \right)_N}   &+ \frac{1}{2}\left\{ {{{\left( {{\nabla  } \cdot \tilde {\mathbf{ F}}\left( {\mathbf{U}} \right),\boldsymbol\phi } \right)}_N} - {{\left( {{\mathbf{U}},{\nabla  } \cdot \tilde {\mathbf{ F}}\left( \boldsymbol\phi  \right)} \right)}_N} + {{\left( {{\nabla  } \cdot {\mathbb{I}^N}\left( {\tilde {\mathcal{A}}} \right){\mathbf{U}},\boldsymbol\phi } \right)}_N}} \right\} \\&+ \int_{\partial E,N} {{{\left\{ {{{\tilde {\mathbf{F}} }^*} - \frac{1}{2}\tilde{\mathbf{ F}} \cdot \hat n} \right\}}^T}\boldsymbol\phi dS}= 0.
\end{split}
\end{equation}

Any of the four equivalent forms $[W] \Leftrightarrow [S] \Leftrightarrow [SC] \Leftrightarrow [DS]$ can be used as is convenient for computation or theory. For instance, to show conservation, choose the form $[W]$ and selectively set each component of $\boldsymbol\phi$ to one. Then the first inner product in the braces vanishes and the second and third cancel leaving
\begin{equation}
\int_{E,N} {J{\mathbf{U}}_{t}d\vec \xi }  =  - \int_{\partial E,N} {{{\contravec {\mathbf{F}} }^*}dS}.
 \end{equation}
 Summing over all elements, the interior face contributions cancel leaving the global conservation statement
 \begin{equation}
 \frac{d}{dt}\sum\limits_{r = 1}^{N_{el}} {\int_{E,N} {{J^r}{{\mathbf{U}}^r}d\vec \xi } }  =  - \sum\limits_{Boundary \atop Faces} {\int_{\partial E,N} {{{\contravec {\mathbf{F}} }^{*,r}}dS} }.
  \end{equation}
  
  To find conditions under which the approximation is constant state preserving, use the form $[S]$ with $\mathbf{U} = \mathbf{c} = \text{const}$ in all elements. The first and third inner products in the braces vanish provided that $\nabla\cdot\left(\mathbb{I}^{N}\left(\contravec{\mathcal A}\right)\right)=0$, and the second is explicitly zero. Consistency of the numerical flux implies that $\tilde{\mathbf{F}}^{*}\left(\mathbf c,\mathbf c;\hat n\right) = \tilde{\mathbf{F}}\cdot \hat n$. Therefore, $\iprodN{J\statevec U_{t},\boldsymbol\phi} = 0$ for all $\boldsymbol\phi\in\mathbb{P}^{N}$, which implies that at each node $nml$ in each element $r$, $d \mathbf{U}^{r}_{nml}/dt=0$.
  
 Finally, the split form approximation is optimally stable in the sense that with the numerical flux (\ref{eq:RiemannFlux}), the norm of the approximate solution satisfies an energy statement like  (\ref{eq:dDimensionWellPosednessLinear_DAK}). We show stability using $[DS]$ and $\boldsymbol\phi=\mathbf U$. With the substitution, the volume terms represented by the first two inner products in the braces immediately cancel. The third inner product in the braces can be bounded
 \begin{equation}
 {\left( {\nabla  \cdot \left({\mathbb{I}^N}\left( {\tilde {\mathcal A}} \right)\right){\mathbf{U}},{\mathbf{U}}} \right)_N} \leqslant \mathop {\max }\limits_E {{{\left\| \frac{{\nabla  \cdot \left({\mathbb{I}^N}\left( {\tilde {\mathcal A}} \right)\right)}}{J} \right\|}_2}}{\left( {J{\mathbf{U}},{\mathbf{U}}} \right)_N} \equiv 2\hat \gamma {\left( {J{\mathbf{U}},{\mathbf{U}}} \right)_N},
 \end{equation}
and under assumptions on the smoothness of $\spacevec{\mathcal{A}}$ \cite{XIE:2013hb} and positivity of the Jacobian \cite{Kopriva:2014yq} the coefficient $\hat \gamma$ will converge spectrally to $\gamma$. (If the divergence of the interpolant vanishes, then $\hat \gamma = 0$.) With the bound on the divergence of the coefficient matrices,
\begin{equation}
\frac{1}{2}\frac{d}{{dt}}\left\| {\mathbf{U}^{r}} \right\|_{J,N}^2 \leqslant  - \int_{\partial E,N} {{{\left\{ {{{\tilde {\mathbf{F}} }^*} - \frac{1}{2}\tilde {\mathbf F}^{r} \cdot \hat n} \right\}}^T}{\mathbf{U}}^{r}dS}  + \frac{1}{2}2\hat \gamma^{r} \left\| {\mathbf{U}}^{r} \right\|_{J,N}^2.
\end{equation}
The change in the total energy is again found by summing over all the elements. If we introduce the integrating factor $\hat \gamma = \mathop {\max }\limits_r \hat\gamma^{r}$, 
\begin{equation}
\begin{split}
\frac{d}{{dt}}\left( {{e^{ - 2\hat \gamma t}}\sum\limits_{r = 1}^{N_{el}} {\inorm{\statevec U^{r}}_{J,N}^2} } \right) 
\leqslant 
-2{e^{ - 2\hat \gamma t}}
&\left\{ \sum\limits_{Boundary\atop Faces} {\int_{\partial E,N} {{{\left( {{\tilde{\statevec F}^*} - \frac{1}{2}\statevec F \cdot \hat n} \right)}^T}\statevec UdS} }  \right.
\\&- \left. \sum\limits_{Interior\atop Faces} {\int_{\partial E,N} {\left( {{\tilde{\statevec F}^{*,T}}\jump{\statevec U} - \frac{1}{2}\jump{{{\left( {\statevec F \cdot \hat n} \right)}^T}\statevec U}} \right)dS} }  \right\}.
\end{split}
\label{eq:IntFaceTerms}
\end{equation}
The interface and boundary terms on the right hand side of (\ref{eq:IntFaceTerms}) are identical to what appeared in the original DGSEM, (\ref{eq:IntFaceTerms1}). Therefore,
the total energy satisfies 
\begin{equation}
\begin{split}
\inorm{\statevec U(T)}_{J,N}^2 &+ \sum\limits_{Boundary\atop Faces} {\int_{0}^{T}\int_{\partial E,N} {{{\statevec U}^T}{\mmatrix A^ + }\statevec UdSdt} }  
\\&
\leqslant {e^{2\hat \gamma T}}\inorm{\statevec U(0)}_{J,N}^2 + \sum\limits_{Boundary \atop Faces} {\int_{0}^{T}\int_{\partial E,N} {{e^{ - 2\hat \gamma (T - t)}}\statevec g ^T\left| {{\mmatrix A^ - }} \right|{{\statevec g} }dSdt} } .
\end{split}
\label{eq:MDStabilityLinear1_DAK}
\end{equation}
As with the continuous solution, (\ref{eq:dDimensionWellPosednessLinearG=0}), if $\nabla\cdot\mathbb{I}^{N}\left(\contravec{\mathcal A}\right) = 0$, $\hat \gamma = 0$, and
\begin{equation}
\begin{split}
\inorm{\statevec U(T)}_{J,N}^2 + \sum\limits_{Boundary\atop Faces} {\int_{0}^{T}\int_{\partial E,N} {{{\statevec U}^T}{\mmatrix A^ + }\statevec UdSdt} }  
\leqslant \inorm{\statevec U(0)}_{J,N}^2 + \sum\limits_{Boundary\atop Faces} {\int_{0}^{T}\int_{\partial E,N} {\statevec g ^T\left| {{\mmatrix A^ - }} \right|{{\statevec g} }dSdt} } .
\end{split}
\label{eq:MDStabilityLinearNoGrowthSplit_DAK}
\end{equation}
Applying the norm equivalence (\ref{eq:NormEquivalence}), we see that the split form approximation is strongly stable for variable coefficient problems and/or curved elements.

Finally, we can use the stability analysis of the split form approximation to write the conditions needed for the original DGSEM to be stable when the coefficients are variable. For simplicity, let us suppose that $\nabla\cdot\spacevec{\mmatrix A} = 0$ and external boundary states $\statevec g=0$ so that the global energy should not increase and instability is not masked by natural growth. Let us also assume that $\nabla\cdot\mathbb{I}^{N}\left(\tilde{\mmatrix A}\right)=0$ so that $\hat \gamma$ also vanishes. Then by (\ref{eq:MDStabilityLinear1_DAK}), the energy of the split form approximation does not grow. With $[DS]\Leftrightarrow [SC]$ and $[DGSEM]\Leftrightarrow [SC]-[C]$, where $[C]$ is the correction term
\begin{equation}\frac{1}{2}{\left( {\left\{ {{\mathbb{I}^N}\left( {\tilde {\mmatrix A}} \right) \cdot \nabla \statevec U + \nabla  \cdot \left( {{\mathbb{I}^N}\left( {\tilde {\mmatrix A}} \right)} \right)\statevec U - \nabla  \cdot \tilde {\statevec F}(\statevec U)} \right\},\boldsymbol \phi } \right)_N},\end{equation}
the elemental energy for the DGSEM satisfies
\begin{equation}\frac{1}{2}\frac{d}{{dt}}\left\| \statevec U \right\|_{J,N}^2 \le  - \int_{\partial E,N} {{{\left\{ {{{\tilde {\statevec F}}^*} - \frac{1}{2}\tilde {\statevec F} \cdot \hat n} \right\}}^T}\statevec UdS}  + \frac{1}{2}\left|{\left( {\left\{ {{\mathbb{I}^N}\left( {\tilde {\mmatrix A}} \right) \cdot \nabla \statevec U - \nabla  \cdot \tilde {\statevec F}(\statevec U)} \right\},\statevec U} \right)_N}\right|.
\label{eq:DGSEM+Aliasing}
\end{equation}
The term in the braces of the volume term is non-zero unless the product rule holds. Therefore, for the DGSEM to be stable when the coefficients are variable, the surface terms (including the dissipation arising from the physical boundaries seen in (\ref{eq:PhyEnergydissipation}) must be sufficiently large to counteract any destabilizing influence of the volume term, which might require trying more dissipative numerical fluxes than the characteristic upwind flux. Practice has shown that at least at low order one can often find numerical fluxes for which the influence of the surface terms is sufficiently dissipative. But (\ref{eq:DGSEM+Aliasing}) shows that the approximation can be unstable if the aliasing growth contribution is larger than the dissipation contribution from the element faces.

\section{Summary}
In this paper, we described a discrete integral spectral calculus for polynomial spectral methods using Legendre-Gauss-Lobatto quadrature. This calculus allowed us to write and analyze discontinuous Galerkin spectral element approximations in a compact notation consistent with the continuous version. In particular, it is possible to easily derive four algebraically equivalent forms of a split form approximation for linear hyperbolic systems. These four equivalent forms can then be used to show global conservation, constant state preservation (when applicable) and, most importantly, strong stability of the split form approximation for variable coefficient problems on curved elements.
\bibliographystyle{plain}

\bibliography{/Users/kopriva/Dropbox/BibTex/dakBib.bib}

\end{document}